\documentclass[a4paper]{amsart} 
\usepackage{amsfonts,amsmath,amssymb,amscd}
\usepackage[all]{xy} 
\usepackage{latexsym}
\usepackage[T1]{fontenc}

\newcommand{\bz}{\mathcal Z}
\newcommand{\bj}{\mathcal J}

\newcommand{\bl}{\mathcal{L}}

\newcommand{\na}{\nabla}

\newcommand{\Bb}{\mathbb{B}}

\newcommand{\Rr}{\mathbb{R}}

\newcommand{\Zz}{\mathbb{Z}}
\newcommand{\Cc}{\mathbb{C}}
\newcommand{\Cp}{\mathbb{P}^1}
\newcommand{\Qq}{\mathbb{Q}}

\newcommand{\ka}{\kappa}
\newcommand{\al}{\alpha}
\newcommand{\ep}{\epsilon}
\newcommand{\si}{\sigma}

\newcommand{\om}{\omega}

\newcommand{\Lam}{\Lambda}

\newcommand{\kah}{\text{K\"{a}hler}}
\newcommand{\del}{\bar \partial}
\newcommand{\siref}{\si_{k,x}}

\newcommand{\beq}{\begin{eqnarray*}}
\newcommand{\eeq}{\end{eqnarray*}}
\newcommand{\bpr}{\begin{preuve}}
\newcommand{\epr}{\end{preuve}}

\newenvironment{preuve}[1][]
{\vskip 2mm\noindent {\it  D\'emonstration#1.   }}{$\Box$ \vskip 2mm}

\newtheorem{defi}{D\'efinition}
\newtheorem{theo}{Th\'eor\`eme}
\newtheorem{lemme}{Lemme}
\newtheorem{prop}{Proposition}

\begin{document}
\title{Hypersurfaces symplectiques r\'eelles et pinceaux de Lefschetz r\'eels}
\author{Damien Gayet}
\address{ Institut Camille Jordan,
Universit\'e Claude Bernard Lyon 1,
43 boulevard du 11 novembre 1918
69622 Villeurbanne cedex
France} 

\email{gayet@math.univ-lyon1.fr}
\maketitle
\centerline{\textbf{Abstract}}
In a compact  symplectic real manifold, i.e supporting 
an antisymplectic involution, we use Donaldson's construction
to build a codimension 2 symplectic submanifold invariant under
the action of the involution.  If the real part of  the manifold is not empty, and if the symplectic form $\om$ is
entire,  then there is an integer $N$ such that for all $k$ big enough, we can find a hypersurface
Poincar\'e dual of $Nk[\omega]$ such that its real part has at least $k^{\frac{n}{2}}$
connected components, up to a constant independant of $k$, and where $2n$ is the dimension of the ambient manifold. 
Finally we extend to our real case Donaldson's construction of Lefschetz pencils. \\
 
\centerline{\textbf{R\'esum\'e}}
Dans le cadre d'une vari\'et\'e symplectique compacte $X^{2n}$  
r\'eelle, c'est-\`a-dire poss\'edant une involution antisymplectique,
nous utilisons la construction
de Donaldson  pour \'etablir l'existence de
sous-vari\'et\'es symplectiques de codimension 2 invariantes par l'involution. 
Si la partie r\'eelle de la vari\'et\'e est non vide,
 et si la forme
symplectique $\omega$ est enti\`ere, alors il existe un entier $N$ tel que 
pour tout degr\'e $k$ assez grand, il existe une hypersurface
Poincar\'e duale \`a $Nk[\omega]$, telle que sa partie r\'eelle poss\`ede au moins  $k^{\frac{n}{2}}$ composantes
connexes, \`a une constante ind\'ependante de $k$ pr\`es, et o\`u $2n$ est la dimension de la vari\'et\'e ambiante. Enfin nous 
\'etendons  au cas r\'eel les r\'esultats de Donaldson 
sur l'existence de pinceaux de Lefschetz.\\

\textsc{Code mati\`ere AMS}: 14P25, 53D05, 14D06,  32Q15.

\subsection*{Introduction}
Soit $(X^{2n},\om)$ une vari\'et\'e symplectique compacte, 
et supposons qu'il existe
 un fibr\'e en droites complexes $L$ de classe de Chern  
$[\om]$, ce qui revient \`a dire que les p\'eriodes de $[\om]$
sont enti\`eres.  
Fixons de plus une structure presque complexe $J$ compatible
avec $\om$.
Dans [Do1], S.K Donaldson montre qu'il existe une  
suite 
de sections $(s_k)_{k}$ de $L^k$ approximativement J-holomorphe (en abr\'eg\'e AH),
c'est-\`a-dire dont le $\del$ est major\'ee par une constante
ind\'ependante de $k$, et dont la d\'eriv\'ee covariante
est minor\'ee aux endroits o\`u $s_k$ s'annule par
$\eta \sqrt k$, o\`u $\eta $ est une constante ind\'ependante de $k$ et strictement positive.
Ainsi, pour $k$ assez
grand, le lieu des z\'eros de $s_k$ est
une vari\'et\'e lisse de codimension r\'eelle 2. De plus, cette sous-vari\'et\'e tend  
\`a devenir de plus en plus $J$-complexe, si bien qu'elle devient symplectique pour $k$ 
assez grand. \\

On suppose  dans cet article que $(X,\om)$
poss\`ede une structure r\'eelle, c'est-\`a-dire une involution 
$c$ v\'erifiant $c^*\om = -\om$. 
Peut-on alors constuire
une suite de sections $(s_k)_{k}$ de $L^k$ du type Donaldson, et 
dont le lieu d'annulation est invariant par $c$ ?
La r\'eponse est positive, et tient dans le 

\begin{theo} Soit $(X,\om,c)$ une vari\'et\'e symplectique r\'eelle.
 Alors il existe une hypersurface symplectique invariante par $c$. Plus pr\'ecis\'ement, il 
 existe une forme symplectique $\tilde \om$ aussi proche de $\om$ qu'on veut, 
 et un entier $N$ tel que pour tout $k$ assez grand, il existe une hypersurface symplectique r\'eelle 
 Poincar\'e duale \`a $Nk[\tilde \om]$.
\end{theo}
\textbf{Topologie de la partie r\'eelle des hypersurfaces.}
Un probl\`eme en g\'eom\'etrie r\'eelle est de conna\^itre la topologie de la partie r\'eelle, c'est-\`a-dire l'ensemble 
des points invariants par $c$, des objets r\'eels construits. 
Nous obtenons le r\'esultat suivant : 
\begin{theo} Il existe une hypersurface symplectique r\'eelle sans partie r\'eelle. Si la partie r\'eelle
de $X$ est non vide, il  existe une forme symplectique $\tilde \om$ 
aussi proche de $\om$ qu'on veut,  un r\'eel $\ep>0$ et un entier $N$ tels que 
pour tout $k$ assez grand, il existe une hypersurface symplectique r\'eelle 
Poincar\'e duale \`a $[Nk\tilde \om]$ dont dont le nombre de composantes
connexes est au moins $\ep k^{\frac{n}{2}}$. Cette borne est \`a une constante pr\`es optimale, dans le sens o\`u 
pour toute suite de sections $(s_{k})_k$ AH et uniform\'ement $\eta$-transverse, il existe une constante $C$
telle que le nombre de composantes connexes de la partie r\'eelle de $s_{k}^{-1}(0)$ 
 est inf\'erieur \`a $Ck^{\frac{n}{2}}$.
 \end{theo}
\textbf{Cas int\'egrable}. Dans le cas o\`u $J$ est int\'egrable, nous avons la proposition suivante :
\begin{prop}
Soit $(X,\om, J, c)$ une vari\'et\'e projective r\'eelle. Alors toutes les hypersurfaces des th\'eor\`emes pr\'ec\'edents
 peuvent \^etre construites complexes. 
\end{prop}
\noindent
\textbf{Pinceaux de Lefschetz r\'eels}. 
Rappelons que si la dimension de $X$ est $2n$, un syst\`eme de coordonn\'ees complexes $(z_{1}, \cdots, z_{n})$ centr\'ees en un point
$x$ est dit \textit{adapt\'e} si la forme symplectique $\omega$ est $(1,1)$ et strictement positive 
au point $x$ pour
la structure complexe induite par ces coordonn\'ees.
\begin{defi}Un pinceau de Lefschetz symplectique associ\'e
\`a une vari\'et\'e symplectique $(X,\om)$ est la donn\'ee de :

(i) Une sous-vari\'et\'e symplectique $N$ de codimension r\'eelle 4.

(ii) Une application surjective $F : X - N \to \Cp $

(iii) Un nombre fini de points $\Delta \subset M-N$
en dehors desquels $F$ est une submersion.

De plus, ces donn\'ees v\'erifient les mod\`eles locaux suivants :

(iv) Pour tout point $p\in N$, il existe une carte adapt\'ee
$(z_1,   \cdots, z_n)$ pour laquelle la sous-vari\'et\'e $N$ 
a pour \'equation locale $\{z_1 = z_2 = 0\}$ et telle que 
$F = z_2/z_1$.

$(v)$ Pour tout point $p\in \Delta$, il existe une carte adapt\'ee
$(z_1, \cdots , z_n)$ dans laquelle $F$ s'\'ecrit
$F (z) = z_1 ^2 + \cdots + z_n^2 + c$. 

Enfin, si $(X,\om)$ est munie d'une structure r\'eelle $c$, 
on dira que le pinceau est r\'eel si 
$$\bar{F} = F \circ c.$$
\end{defi}
Le th\'eor\`eme principal de [Do2] est le suivant :
\begin{theo}([Do2]) Soit $(X,\om)$ une vari\'et\'e symplectique compacte telle 
que la classe de cohomologie $[\om]$ soit enti\`ere. Alors pour $k$ assez grand,
il existe un pinceau de Lefschetz dont les fibres sont Poincar\'e duales 
de $k[\om]$. 
\end{theo}
Nous  d\'emontrons dans la deuxi\`eme partie de cet article que  ce th\'eor\`eme s'adapte 
dans le cas r\'eel : 
\begin{theo}
Soit $(X,\om,c)$ une vari\'et\'e symplectique r\'eelle, telle que $\om$ soit
enti\`ere. Alors il existe un pinceau de Lefschetz r\'eel.
\end{theo}

\noindent
\textit{Remerciements}: je remercie Jean-Yves Welschinger pour 
m'avoir pos\'e cette question et pour ses judicieux commentaires tout au long de ce travail, ainsi 
que le referee pour ses suggestions et corrections.

\section {Construction d'une hypersurface symplectique r\'eelle}
\subsection {Structures associ\'ees \`a la structure r\'eelle}
\noindent
\textbf{Structures presque complexes. }Une premi\`ere remarque d\'emontr\'ee dans [We] est 
qu'il existe une structure presque complexe $J$
\`a la fois compatible avec $\om$ 
et rendant $c$ antiholomorphe, c'est-\`a-dire
telle que $dc\circ J = - J dc$. Dans toute la suite, $J$ d\'esignera une telle structure.\\
\noindent
\textbf{Forme symplectiques \`a valeurs enti\`eres. }
Pour que $-i\om$ soit la courbure d'un fibr\'e en droites complexes, 
il faut que ses valeurs sur les 2-cycles soient enti\`eres. Le lemme
suivant nous permet de r\'ealiser cette condition, tout en restant dans un 
cadre r\'eel. 
\begin{lemme} 
Si $\om $ n'est pas \`a valeurs rationnelles, 
on peut la perturber en $\tilde \omega$ de sorte
que $c$ reste antisymplectique pour la nouvelle forme.
\end{lemme}
\bpr
En effet, soit $\bz_{+}^2$ le sous-espace vectoriel 
propre pour la valeur propre $+1$ de l'endomrophisme $-c^*$:
 $$-c^*: \bz^2(X,\Rr) \to \bz^2(X,\Rr)$$ 
agissant sur les 2-formes ferm\'ees r\'eelles. 
Le r\'eseau $\bz^2(X,\Qq)$ des 2-formes ferm\'ees
\`a valeurs rationnelles sur les 2-cycles est invariant par $-c^*$, et son intersection
$\bz^2_{+}\cap \bz^2(X,\Qq)$  est dense dans $\bz^2_{+}$. On peut 
donc perturber $\om$ en un \'el\'ement de $\bz^2(X,\Qq)$ tout en restant
dans l'espace des formes invariantes par $-c^*$. 
\epr

\noindent
\textbf{Fibr\'es  $c$-r\'eels.}
 Le lemme suivant est imm\'ediat : 
\begin{lemme}
Soit $E$ un fibr\'e  complexe au-dessus de $X$. Alors l'application $c$
de $X$ dans $X$ se rel\`eve en un isomorphisme 
$\Cc$-antilin\'eaire de fibr\'es $\hat{c} : E \to E$ 
si et seulement si $(c^*E)^* = E$.
\end{lemme}
Le probl\`eme est qu'il n'existe pas toujours une $involution$ $\Cc$-antilin\'eaire qui 
rel\`eve $c$, voir par exemple [We2]. 
\begin{defi} Un fibr\'e $E$ est dit $c$-r\'eel s'il existe une
\textit{involution} antilin\'eaire de $E$ dans $E$ relevant $c$. Dans ce cas, 
on d\'esigne par $\kappa : \Gamma(E) \to \Gamma (E) $ l'involution sur l'espace des sections
de $E$ induite par $\hat{c}$ : $$\kappa (s) = \hat{c}^{-1}\circ s \circ c.$$
Une section $s$ de $E$ est dite
sym\'etrique si $\ka(s)=s$.
\end{defi}
\begin{lemme}
Soit $L$ un fibr\'e en droites v\'erifiant $L=(c^*L)^*$. Alors $L^2$ est $c$-r\'eel.
\end{lemme}
\begin{preuve}
Soit $\hat{c}$ une application  antilin\'eaire de $L$ dans $L$ relevant $c$.
 L'application $\hat{c}\circ \hat {c}$ est un isomorphisme lin\'eaire de $L$ dans $L$ relevant
l'identit\'e. Il existe donc une fonction complexe $a: X \to \Cc^*$ 
telle que  $\hat{c}\circ \hat{c}(x,\lambda) = (x,a(x)\lambda)$. 
En utilisant l'\'egalit\'e $\hat{c}^2 \circ \hat{c} = \hat{c} \circ \hat{c} ^2$, on constate que 
cette fonction $a$ v\'erifie l'identit\'e : $$a = \overline{a(c)}.$$
Si l'on veut construire un nouveau rel\`evement de $c$, il suffit de diviser $\hat{c}$ par une application 
complexe ne s'annulant pas : $\hat{c}_{\phi} = \frac{1}{\phi}\hat{c}$. Ainsi, 
$\hat{c}_{\phi}^2 = \hat{c}^2 \frac{1}{\phi \overline{\phi (c)}}.$
Maintenant, il est facile de v\'erifier que l'application 
$$\hat{c}_{L^2} = \frac{1}{a}\ \hat{c}\otimes \hat {c}$$
est une involution antilin\'eaire de $L^2$ dans $L^2$ relevant $c$.
\end{preuve}
\noindent
\textbf {Fait :} Soit $(X^{2n},J)$ une vari\'et\'e presque complexe munie
d'une involution $J$-antiholomorphe. Alors pour tout $p$ et $q$ inf\'erieurs \`a $n$, 
le fibr\'e des $(p,q)-$formes est $c$-r\'eel, et ont peut choisir $\kappa(\beta) = \overline{c^*\beta}$\\

Supposons que la partie r\'eelle  $\bl = \Rr X$ soit non vide. Dans ce cas,
il est facile de d\'emontrer que  $\bl$
est une sous-vari\'et\'e lagrangienne. Par cons\'equent, la restriction 
du  fibr\'e en droites complexes $L$ \`a $\bl$ est topologiquement triviale
puisque sa classe de Chern $[\om_{|\bl}]$ est nulle.   Quitte \`a remplacer
$L$ par une certaine puissance $L^N$, il existe donc une section $e$ de $L$
ne s'annulant pas sur un  voisinage $V$ de $\bl$. Ceci nous permet
de choisir une  involution de $L$ particuli\`ere :
\begin{lemme}
Si $\Rr X$ est non vide, et quitte \`a remplacer $L$ par une puissance paire $L^{2N}$,
 pour toute section 
trivialisante $e$ de $L$ sur le voisinage $V$ de $\Rr X$, de norme 1
sur $V$, on peut choisir $\hat{c}$ de sorte que $\hat{c} (e) = e$ aux points 
de $\Rr X$. 
\end{lemme}
\bpr
Soit $f \in C^\infty(V,S^1)$, telle que
$\hat{c} (e) = f e $ sur $V$. 
Etendons 
\`a $X$ cette application. Quitte \`a consid\'erer
$L^2$ au lieu de $L$, l'application $ff(c)$ admet une racine
carr\'ee $g$ qu'on peut choisir \'egale \`a $f$ sur $\Rr X$ et v\'erifiant $g(c)=g$. Maintenant, 
le nouveau morphisme de fibr\'e $\tilde c = \frac{1}{g} \hat c$ est l'identit\'e au-dessus
de $\Rr X$. 
\epr

\subsubsection{Connexions r\'eelles}
On suppose dor\'enavant que $L$ est un fibr\'e $c-$r\'eel. Sur $L$
on fixe une m\'etrique hermitienne invariante par $c$, qui existe toujours. 
L'involution $\kappa$ s'\'etend \`a  $\Gamma(TX^*\otimes L)$ par :
$$\kappa : \alpha\otimes \lambda \mapsto \overline{c^{*}\alpha}\otimes \kappa (\lambda).$$
Si $\nabla $ est une connection hermitienne, $\nabla s \in \Gamma (TX^*\otimes L)$, 
ce qui  nous permet de construire une involution naturelle sur les
 connexions  de $L$ par 
 $$\ka (\nabla) =\kappa \nabla \kappa.$$Il est alors ais\'e de d\'emontrer le lemme suivant : 
\begin{lemme} Si $\nabla$ est une connection hermitienne de courbure $-i\om$, alors $\kappa (\nabla)$ \'egalement. 
De plus pour toute connexion $\ka$-invariante,
et toute section $s$ de $L$, on a $|\nabla \kappa s | = |\nabla s| (c)$, de m\^eme que 
$|\na^{0,1}\ka (s)| =|\na^{0,1}s|(c).$ 
\end{lemme}
On a par ailleurs le 
\begin{lemme} Il existe une connexion unitaire $\na$ de courbure $-i\omega$ et $\ka$-invariante 
sur $X$.
\end{lemme}
\bpr Si $\nabla_{0}$ est une connexion hermitienne de courbure $-i\omega$, alors 
$\nabla= \frac{1}{2}(\nabla_{0}+ \ka(\nabla_{0}))$ convient. 
\end{preuve}
Dans toute la suite, on munira notre fibr\'e d'une telle connexion.

\subsection {Le th\'eor\`eme de Donaldson [Do1]}
Soit $(X,\om, J)$ une vari\'et\'e symplectique compacte munie d'une structure
presque complexe $J$ compatible avec $\om$, c'est-\`a-dire telle que
$\om(.,J.)$ est une m\'etrique riemanienne. Si $\om$ est \`a p\'eriodes
enti\`eres, il existe un fibr\'e en droites complexes  $L$ sur $X$, 
dont la classe de Chern repr\'esente la classe $[\om]\in H^2(X,\Zz)$. 
Si l'on muni $L$ d'une norme hermitienne, il existe une connexion
unitaire  $\nabla$ sur $L$ de courbure $-i\om$. Tous ces objets
ont une extension naturelle pour chaque puissance $L^k$ du fibr\'e.
En particulier, on notera $g_{k}$ la m\'etrique $kg$. 

\begin{defi} Soit $E$ un fibr\'e hermitien sur $X$. 
Une suite $(s_k)_{k}$ de sections de $E\otimes L^k$ 
est asymptotiquement holomorphe (AH)
s'il existe une constante $C$ ind\'ependante de $k$, telle que pour
tout $k$, $|\nabla^{0,1} s_k|\leq C$, telle que
 sa d\'eriv\'ee est en norme inf\'erieure \`a $C\sqrt k$, et sa d\'eriv\'ee seconde en est de norme
 inf\'erieure \`a $Ck$. La suite $(s_k)_k$
est de plus $\ep$-transverse sur $X$ s'il existe $\ep>0$, tel que pour
tout $x$ tel que $|s_k(x)|\leq \ep$, $\nabla s_k(x)$ 
admet un inverse \`a droite de norme inf\'erieure \`a $(\eta\sqrt k)^{-1}$. 
\end{defi}
Le th\'eor\`eme principal de [Do1], originellement d\'emontr\'e si $E=X\times \Cc$, 
et g\'en\'eralis\'e par Auroux pour le cas g\'en\'eral, est le suivant : 
\begin{theo}(Donaldson, Auroux) Soit $E$ un fibr\'e hermitien sur $X$. 
Pour toute suite de sections $(s_k)_{k}$ AH de $E\otimes L^k$, et pour tout $\ep>0$ fix\'e,
alors il existe un $\eta>0$, une suite de sections $(\tilde s_k)_{k}$ AH, $\eta$-transverse
pour $k$ assez grand, 
avec $|s_k - \tilde s_k|\leq \ep$, et $|\nabla s_k - \nabla \tilde s_k|\leq \ep \sqrt k$. 
\end{theo}

Il nous faut rappeler les id\'ees de la d\'emonstration de ce th\'eor\`eme, dans le cas o\`u $E$ 
est le fibr\'e en droites complexes trivial. 
Pour tout r\'eel $D>0$ a priori fix\'e, il existe un recouvrement 
par ``couleurs'', c'est-\`a-dire par des ensembles
de boules de $g_k$-taille \'egale \`a 1 et
 dont les centres sont $g_k$-distants d'au moins $D$. 
Le nombre $N$ de ces couleurs augmente \'evidemment avec $D$, mais pas avec $k$. On perturbe 
en $N$ \'etapes la suite de sections $s_k$ par une somme 
pond\'er\'ee de sections essentiellement localis\'ees au-dessus des centres 
des boules d'une couleur donn\'ee. Plus pr\'ecis\'ement, ces sections sont du m\^eme mod\`ele :
\begin{lemme}
Soit $x\in X$. Il existe une suite de sections $(\si_{k,x})_{k}$ asymptotiquement
holomorphe satisfaisant les in\'egalit\'es suivantes :
\begin{enumerate}
\item   $|\siref| \geq \frac{1}{2} sur \  B_{g_k}(x,1)$
\item  $|\siref|_{C^2}\leq p(d_k(x,y))e^{-d_k(x,y)^2}$
\item  $\siref$ \text{ a  un support  inclus  dans } $B_{g_k}(x, k^{1/6}) = B_{g}(x, k^{-1/3})$,
\end{enumerate}
 o\`u $p$ est un polyn\^ome ind\'ependant de $k$. 
\end{lemme}
La perturbation sur une couleur est faite de la façon suivante :
$$ s_k \rightarrow s_k+\sum_{x_i\in couleur}w_i \si_{k,x_i}.$$
A chaque \'etape, un th\'eor\`eme de Sard permet de choisir les $w_{i}$ de façon \`a ce que 
la nouvelle section devienne $\eta$-transverse sur la nouvelle
couleur et le reste sur les anciennes, mais pour un $\eta$ de plus en plus 
petit. Le th\'eor\`eme de Sard utilis\'e est quantitatif, ce qui permet d'effectuer cette r\'ecurrence 
en un nombre d'\'etapes ne d\'ependant pas de $k$. 

\subsection{Construction de l'hypersurface r\'eelle}
La plupart des outils techniques qui nous permettent
d'obtenir le th\'eor\`eme 1 sont  partiellement pr\'esents dans
les articles [Do1], [Do2], [Au] et [Au,Mu,Pr]. Dans ce dernier article il est d\'emontr\'e
que si $\bl$ est une sous-vari\'et\'e lagrangienne
munie d'une fonction de Morse r\'eelle, alors 
il existe un pinceau de Lefschetz $F$ \`a la Donaldson,
tel qu'une isotopie de  $\bl$ se projette par $F$ sur
un arc r\'eel. Pour d\'emontrer ce r\'esultat, les auteurs ont 
besoin de consid\'erer $\bl$ comme la partie
r\'eelle d'une stucture r\'eelle semi-locale, et de construire
les perturbations de façon sym\'etriques par rapport \`a cette structure au voisinage
de $\bl$. Dans notre cas, nous devons construire 
les sections et leur perturbation de façon sym\'etrique, mais cette fois globalement sur tout $X$.
Par ailleurs, nous aimerions transversaliser n'importe quelle suite de sections AH sym\'etrique, 
pas n\'ecessairement \`a partir d'une fonction de Morse sur $\Rr X$. 
Nous allons en fait d\'emontrer la proposition
suivante.
\begin{prop} Soit $E$ un fibr\'e hermitien $c-$r\'eel, et 
 $(s_k)_{k}$ une suite de sections sym\'etriques et AH de $E\otimes L^k$.  Alors il est possible 
de la perturber en une suite de sections sym\'etriques, AH et 
transverse sur $X$.
\end{prop}
Remarquons que le th\'eor\`eme 1 est une cons\'equence imm\'ediate
de cette proposition, avec $E=X\times \Cc$.\\
\textbf{Fait : } on peut toujours choisir le r\'eseau de boules, ainsi que chaque couleur, 
invariante par $c$. De plus, on peut faire en sorte que si l'une de ces boules rencontre $\Rr X$, alors
la boule est invariante par $c$. 


\noindent
Par souci de clart\'e, nous d\'emontrons le th\'eor\`eme pour $E=X\times \Cc$, puis nous 
expliquerons les adaptations n\'ecessaires pour le cas g\'en\'eral. 

\noindent
\bpr[ dans la cas o\`u $E=X\times \Cc$] 
D\'esignons par $\Lam'$ (resp. $I'_i$) un sous-r\'eseau 
minimal  (resp. de la couleur $I_i$)
tel que $$\bigcup_{x \in \Lam'} B_{g_k}(x, 1) \cup B_{g_k}(c(x)), 1) = X,$$
(resp. pour la couleur $i$, \textit{mutatis mutandis}). 
Lors de la transversalisation sur une boule de la couleur $I_i$, 
si l'on perturbe  la section $s_k$ par $w \si_{k,p}$,
on perd le caract\`ere r\'eel. Nous avons donc besoin d'un raffinement 
du  th\'eor\`eme de Sard quantitatif pr\'esent dans [Do1] :
\begin{prop} 
Il existe un entier $p$ et $\delta_0>0$,  tels que pour tout $0<\delta<\delta_0$, 
si $\si = \delta \log (\delta^{-1})^{-p}$,  et  $f$ une 
fonction complexe  d\'efinie sur la boule
 $B(0,\frac{11}{10})  \subset \Cc^n$ v\'erifiant
$$ |f|\leq 1, \ \ |\del f|_{C^1}\leq \si,$$
alors il existe une  constante $w\in \Rr$ , avec $|w|\leq \delta$, et
$f-w$ est $\si$-transverse sur la boule unit\'e de $\Cc^n$. 
\end{prop}
\noindent
\textbf{Remarque :} Notons que la seule diff\'erence avec les propositions de [Do1]
est que l'on peut choisir $w$ r\'eel.\\ 
Pour poursuivre, nous avons besoin \'egalement du raffinement  suivant du lemme 7, qui nous
permet de remplacer $\si_{k,x}$ par une section AH sym\'etrique  : \\
\textbf{Fait : } Pour tout $x\in \Lam'$, la section AH et sym\'etrique 
$$\hat \si_{k,x} = \frac{\si_{k,x}+ \ka( \si_{k,x})}{2}$$
est de norme  uniform\'ement minor\'ee sur $B_{g_{k}}(x,1)\cup B_{g_{k}}(c(x),1)$.\\

On applique donc le proc\'ed\'e de Donaldson sur les boules centr\'ees sur les points $x$
 de $\Lam'$, mais en utilisant la proposition 3 appliqu\'ee \`a 
$$f = \frac{s_k}{\hat\si_{k,x}},$$
ce qui permet de trouver un $w$ $r\'eel$ de sorte que  
  $ s_k + w \hat \si_{k,x}$ est sym\'etrique et transverse sur la boule $B_{g_k}(x,1)$. 
 On a alors par le lemme 5 la transversalisation
automatique sur les boules centr\'ees sur les points de $c(\Lam')$, ce qui d\'emontre le th\'eor\`eme 1. 
\epr

Nous donnons maintenant une d\'emonstration du raffinement du lemme de Sard.
\bpr [ de la proposition 3] 
La preuve suit jusqu'au dernier moment  la preuve classique donn\'ee par [Do1]. L'application
$f$ v\'erifiant  $|\del f| \leq \si$ peut \^etre $C^1$-approch\'ee par une application holomorphe
$\tilde f$ sur la boule unit\'e, et l'erreur est inf\'erieure \`a un multiple (born\'e) de $\si$. 
Ensuite
il est possible de $C^1$-approcher $\tilde f$ de $\si$ par un polyn\^ome complexe 
$g$ de degr\'e inf\'erieur ou 
\'egal \`a  $C\log (\delta^{-1})$, o\`u $C$ est une constante ne d\'ependant que de la 
boule. Pour toute fonction complexe $h$, 
soit $$Y_{h,\epsilon}=\{z\in \Bb^{2n}, \ |d_{x}h|\leq \epsilon \},$$
et $Z_{h,\epsilon}$ le $\ep$-voisinage tubulaire de $f(Y_{h,\epsilon})$. 
Nous avons l'inclusion : $Z_{f,\sigma }\subset Z_{g,c\sigma}$, o\`u $c$ est une
constante de structure. 
Si l'on trouve un r\'eel $w$ dans un disque de taille $\delta$
\'evitant $Z(g,c\sigma)$, alors on a prouv\'e la proposition.
Rappelons la 
\begin{prop}Soit $P$ un polyn\^ome de $\Rr^m$ dans $\Rr$, 
tel que 1 soit une valeur r\'eguli\`ere de P sur $\Bb^m$, ainsi
que de la restriction de P sur $S^{m-1}$. Alors 
le nombre de composantes du sous-niveau $\{P\leq 1\}$ ainsi
que leur diam\`etre (pour la m\'etrique induite) sont $p-$born\'es.
\end{prop}
On dit qu'une quantit\'e associ\'ee \`a un polyn\^ome est $p$-born\'ee si elle est major\'ee
par une puissance du degr\'e du polyn\^ome. 
Apr\`es une infime perturbation de $g$, appliquons la proposition
pr\'ec\'edente au polyn\^ome r\'eel $|\partial g/c\sigma|^2$, qui est de 
degr\'e $C\log (\sigma^{-1})$, o\`u $C$ est une nouvelle constante de structure. 
Le diam\`etre pour la m\'etrique ambiante de l'image d'une composante 
de $Y_{g,c\sigma}$ est major\'e par $C'\sigma \log (\sigma^{-1})^{p}$, o\`u $p$
est la puissance donn\'ee par la proposition. 
Le voisinage $\si$-tubulaire de $Z_{g,c\sigma}$ est donc contenu dans 
une r\'eunion de $\log (\sigma^{-1})^{p}$
 disques de rayons \'egaux \`a  $\sigma \log (\sigma^{-1})^{p}+\si$. 
 Ce dernier ajout de $\sigma$ ne change rien, et on peut l'oublier. 
L'intersection de l'ensemble prohib\'e 
avec $\Rr$ recouvre donc au plus une longueur
de  $\sigma \log (\sigma^{-1})^{2p}$. 
Si le rayon du disque est par exemple $\sigma \log (\sigma^{-1})^{3p}$,
il reste encore beaucoup de place pour trouver un $w$ r\'eel 
dans le compl\'ementaire de $Z_{g,c\sigma}$ dans $[-\delta, \delta]$. 
On obtient alors la proposition, en effectuant un changement de variables
comme dans [Do1]. \\

\noindent
\textit{D\'emonstration de la proposition 2 dans le cas o\`u $E$ est
un fibr\'e hermitien quelconque.}
Nous suivons la d\'emonstration de [Au] et l'adaptons en cours de route \`a notre cas. Elle consiste \`a 
transversaliser $s_k$ composante par composante. 
Plus pr\'ecis\'ement,
on choisit un recouvrement par des ouverts $U_i$
de trivialisation du fibr\'e $E$,   tels que $c(U_i)$ coïncide avec
un autre $U_j$, et pr\'ecis\'ement $U_i$ si $U_i$ rencontre $\Rr X$. 
On transversalise sur la moiti\'e de l'ensemble de ces paires d'ouverts en ajoutant 
\`a $s_k$ des sections sym\'etriques, si bien qu'automatiquement,
cela donnera la transversalisation sur le reste des ouverts. Supposons que sur un certain $U_i$, les $r$ premi\`eres 
composantes de $s_k$, qu'on \'ecrit $\pi_{\leq r} s_k$,
 forment une section transverse, et 
d\'efinissent donc une sous-vari\'et\'e symplectique $W_{k,r}$
lisse sur $U_i$.
 Auroux montre qu'
il est possible de transversaliser  la restriction de la $r+1$-i\`eme composante
$\pi_{r+1} s_k$ sur $W_{k,r}$ en lui ajoutant
une petite section AH $\tau_{k,r}$. Il d\'emontre ensuite que cela implique automatiquement
que $\pi_{\leq r}s_k\oplus (\pi_{r+1} s_k + \tau_{r,k})$ est en fait
transverse sur tout $U_i$. Dans notre cas, nous savons que $W_{k,r}$ est invariant
par $c$. 
 Par ailleurs, puisque la section partielle $\pi_{\leq r} s_k$
est transverse, il est facile de voir qu'il existe un $\rho>0$ ind\'ependant de $k$, 
tel que pour tout point  $x\in \Rr X$,  l'intersection de $B_{g_k}(x,\rho)$
avec $W_{k,r}$ est toujours ou bien  connexe et de partie r\'eelle non vide, ou bien est vide. 
Pour le proc\'ed\'e de Donaldson, on consid\`ere donc un r\'eseau $\Lambda'$ invariant par $c$, et centr\'e
sur des points de $W_{k,r}$,  de boules $B_{g_k}(x,\rho)$.
 Il est alors possible d'appliquer la d\'emonstration 
d\'ecrite dans [Au], sans param\`etre $t$, et en utilisant la proposition 3
pour chaque transversalisation, et enfin  $\hat \si_{k,x}$ au lieu de $\si_{k,x}$.  
 Au total, on a ajout\'e
des sections sym\'etriques qu'on somme en une section $\tau_{k,r}$, 
de sorte que la section 
$\pi_{\leq r}s_k\oplus (\pi_{r+1} s_k + \tau_{r,k})$ est sym\'etrique et
transverse sur $U_i$. Pour tout ouvert $U_i$ trivialisant
$E$, on r\'ep\`ete le processus, dont le nombre de pas ne d\'epend que de $E$ et $X$, si bien 
qu'au total, on obtient la proposition 2 dans toute sa g\'en\'eralit\'e.
\epr
\subsection{Unicit\'e}
Nous d\'emontrons maintenant une version d'unicit\'e des hypersurfaces r\'eelles construites : 
\begin{prop}
Soient $(s_{1,k})_k$ et $(s_{2,k})_k$ deux suites de sections  de $E\otimes L^k$ AH, transverses et sym\'etriques, et 
telles que $s_{1,k}(x)=s_{2,k}(x)$ pour tout $x\in \Rr X$. 
Alors pour tout $k$ assez grand, il existe une isotopie de sections AH et transverses reliant 
les lieux d'annulations de $s_{1,k}$ et $s_{2,k}$. 
\end{prop}
\bpr
Pour faciliter la lecture, nous d\'emontrons cette proposition dans le cas o\`u $E=X\times \Cc$. 
Tout comme dans [Au], on consid\`ere la suite de sections sym\'etriques
$s_{t,k}= ts_{1,k}+(1-t)s_{2,k}$,
et l'on va montrer qu'on peut transversaliser cette suite 
uniform\'ement en $t$. Le probl\`eme principal est que la proposition 3, contrairement
\`a sa version plus souple de [Do1], n'admet pas 
de g\'en\'eralisation avec un param\`etre. En effet, 
le lieu des $w$ r\'eels interdits \`a chaque temps $t$ 
est certes tr\`es petit, mais s\'epare en g\'en\'eral l'axe des r\'eels
en au moins deux parties, au sein desquelles il faut choisir
$w_t$. Or l'une de ces parties peut dispara\^itre pour un certain
$t$, si bien qu'il faudrait alors choisir des $w_t$ discontinus, 
ce qui n'est \'evidemment pas possible. 
 Nous verrons plus bas que la condition d'\'egalit\'e des deux sections aux points de $\Rr X$ permet 
de s'abstenir de cette version. En dehors d'un voisinage de $g$-taille fixe de $\Rr X$, 
nous pouvons utiliser la version param\'etrique d'Auroux (proposition 3 dans [Au]), 
et additionner \`a $s_{t,k}$ la section $w_t\si_{k,x}+ \overline{w_t} \ka(\si_{k,x})$. 
Cette derni\`ere est sym\'etrique, et puisque les supports de $\si_{k,x}$ et $\ka (\si_{k,x})$ 
sont disjoints, d'une part  la transversalit\'e aquise d'un c\^ot\'e, sur $B_{g_k}(x,1)$,
n'est pas perturb\'ee par le terme $\overline{w_t}\ka(\si_{k,x})$, d'autre part par le lemme 5,  
la transversalit\'e est spontan\'ement acquise de l'autre c\^ot\'e, sur $c(B_{g_k}(x,1))$ par ce dernier
terme. 

Maintenant, \'etudions la situation sur une boule $B_{g_k}(x,1)$, avec $x\in \Rr X$. 
Apr\`es trivialisation par le lemme 9 ci-dessous, nous utilisons le lemme suivant : 
\begin{lemme}
Soit $f$ une fonction complexe d\'efinie sur $2\Bb^{2n}$,  v\'erifiant
$|\del f|_{C^1} < \ep$, et $f (x) = 0$ pour tout $x$ dans $\Rr^n$.
Alors $|f|_{C^1(\Bb^{2n})} < c\ep$, o\`u $c$ est une constante ind\'ependante 
de $f$ et de $\ep$. 
\end{lemme}
\bpr Pour d\'emontrer ce lemme, on commence par trouver gr\^ace
au lemme 28 de [Do1] une fonction holomorphe $\tilde f$
sur $\frac{3}{2}\Bb^{2n}$ telle que $${|f-\tilde f|_{C^1(\frac{3}{2}\Bb^{2n})}< K\ep},$$ o\`u $K$ est une constante
ind\'ependante de $f$ et de $\ep$. 
Ensuite, le lemme 27 de [Do] nous donne un polyn\^ome complexe $p$, 
tel que $|\tilde f - p|_{C^1}<\ep$. On a donc par la seconde hypoth\`ese 
$|p_{|\Rr^n}|_{C^1}<(K+1)\ep$, 
ce qui implique, puisque $p$ est complexe, 
l'existence d'une constante $c>0$ telle que  $|p|_{C^1}<c(K+1)\ep$,
et donc $|f|_{C^1}<(1+c)(K+1)\ep$. 
\epr
En appliquant ce lemme \`a $f=\frac{s_{1,k}- s_{2,k}}{\si_{k,x}}$ via la traditionnelle
renormalisation
par $1/\sqrt k$, 
on obtient que $s_{1,k}- s_{2,k} $ est en norme inf\'erieure \`a $C/\sqrt k$, 
et sa d\'eriv\'ee covariante inf\'erieure en norme \`a $C$ (pour la m\'etrique $g$). 
C'est donc \'egalement vrai pour  $s_{t,k} -s_{2,k}$
uniform\'ement en $t$, et puisque $s_{2,k}$
est transverse sur $B_{g_k}(x,1)$, la section
$s_{t,k}$ l'est aussi, pour $k$ assez grand. 

Il reste maintenant \`a transversaliser $s_{t,k}$ 
dans l'espace entre les deux r\'egions pr\'ec\'edentes. Pour cela, il nous faut un raffinement 
de la proposition 3 de [Au], qui est donn\'e par la proposition suivante de [Au,Mu,Pr] :
\begin{prop}{[Au,Mu,Pr]}
Soit $C>0$, $\ep>0$. Alors il existe un entier $p$ et une constante
$\delta_0$, tels que pour tout $0<\delta<\delta_0$, si 
$\si = \delta (\log(\delta^{-1}))^{-p}$ et $f_t$, $h_t$ deux
 chemins continus 
 param\'etr\'es par $[0,1]$ de fonctions complexes sur la boule $B(0,\frac{11}{10})\subset \Cc^n$
v\'erifiant les majorations suivantes :
$$
|f_t|\leq 1, \ \  |\del f_t|_{C^1}\leq \si,  \ \
|h_t| \leq 1-\ep, \ \ |dh_t|\leq C, \ \ |\del h_t|_{C^1}\leq \si.
$$
Alors il existe un chemin continu $w_t\in \Cc$ v\'erifiant $|w_t|< \delta$, tel que $f_t-w_t-\bar w_t h_t$
soit $\si$-transverse \`a 0 sur la boule unit\'e de $\Cc^n$. 
\end{prop}
Maintenant,  il suffit d'ajouter \`a $s_{t,k}$ la section 
$w_t\si_{k,x} + \bar w_t\ka(\si_{k,x})$, o\`u $w_t$
est d\'etermin\'ee par la proposition pr\'ec\'edente. 
En effet, la section 
$s_{t,k} -w_t\si_{k,x} - \bar w_t \kappa (\si_{k,x})$ correspond \`a
la fonction  $$f_t = \frac{ s_{t,k}}{\si_{k,x}} - w_t -\bar w_t \frac{\ka(\si_{k,x})}{\si_{k,x}}.$$
Quitte \`a prendre des boules de trivialisations assez petites (de $g_k$-taille ind\'ependantes de 
la boule et de $k$), les deux points $x$ et $c(x)$ sont \`a une $g_k$-distance
d'au moins $2$, si bien que la fonction $h_t = h = \frac{\ka(\si_{k,x})}{\si_{k,x}}$ et $f_t$
v\'erifient les hypoth\`eses de la proposition. 
\epr

\subsection{La partie r\'eelle de l'hypersurface symplectique}
Dans ce paragraphe, on suppose que la partie r\'eelle est non vide. Le but
est de construire des  hypersurfaces r\'eelles symplectiques
dont le lieu r\'eel est ou bien vide, ou bien non vide, et 
dans ce cas de maximiser si possible sa topologie. 
Nous pouvons \'enoncer deux propositions simples.
Le premier cas est minimal, puisque la partie r\'eelle est 
vide : 
\begin{prop}
Il existe une hypersurface symplectique de partie r\'eelle vide.
\end{prop}
\bpr
Dans [Au,Ga,Mo], les auteurs ont construit une section AH de norme uniform\'ement minor\'ee au-dessus d'une
sous-vari\'et\'e lagrangienne $\bl $ donn\'ee. Ce th\'eor\`eme est raffin\'e dans le lemme 5.4 de [Au,Mu,Pr], o\`u la section est cette fois
localement sym\'etrique au voisinage de $\bl$ si celle-ci est localement la partie r\'eelle d'un involution antisymplectique, et s'annule
au-del\`a de ce voisinage, ce qui convient \`a notre situation. Il suffit alors de transversaliser par notre proposition 2 pour
obtenir une section sym\'etrique et ne s'annulant pas sur $\Rr X$. 
\epr
La proposition suivante cherche, au contraire, \`a maximiser la topologie de la partie r\'eelle de l'hypersurface :  
\begin{theo}
Soit $(X^{2n},\om, J, c)$ une vari\'et\'e symplectique r\'eelle de partie r\'eelle non vide.  Si $[\om]$ 
est rationnelle, il existe un $\ep >0 $ tel que pour tout $k$ assez grand, 
il existe une hypersurface symplectique r\'eelle Poincar\'e duale \`a $2k[\om]$ et une suite de sections AH $(s_{k})_{k}$, tel que le 
nombre de composantes connexes de sa partie r\'eelle est  au moins $\ep  k^{\frac{n}{2}}$.  
\end{theo}
\noindent
\textbf{Remarque 1.} Rappelons que dans [Do1], il est d\'emontr\'e que
les hypersurfaces sont toujours connexes pour $n\geq 2$. \\

\noindent 
\textbf{Remarque 2.} En dimension 4,  notre th\'eor\`eme donne une minoration en $\ep k$, ce qui est \`a comparer avec 
le th\'eor\`eme d'Harnack, qui majore  dans $\Cc P^2$ complexe 
le nombre d'ovales d'une courbe holomorphe r\'eelle par une borne de l'ordre de $k^2/2$.  \\

\noindent
\textbf{Remarque 3.} Les hypersurfaces de Donaldson tendent \`a remplir tout l'espace, et donc n'ont pas de raison
particuli\`ere de se concentrer sur la partie r\'eelle de $X$. On a en effet ( [Do1], proposition 40) :
$$ \frac{1}{k}Z_{k} \to \om$$
en tant que courants. On peut donc s'attendre, dans notre cas,  \`a ce que la partie r\'eelle de ces hypersurfaces
soit relativement repr\'esentative de celle du lieu d'annulation d'une section prise au hasard. \\

\noindent
\textbf{Remarque 4. }Dans [Ed,Ko], les
auteurs d\'emontrent que dans $\Cc P^1$ et pour une certaine mesure naturelle sur les polyn\^omes de degr\'e $k$,
le nombre moyen de racines r\'eelles d'un polyn\^ome r\'eel est $\sqrt k$, ce qui est pr\'ecis\'ement 
notre cas \`a une constante pr\`es. \\
 \noindent
 La proposition suivante montre que de toutes façons il ne fallait pas s'attendre \`a obtenir vraiment mieux que 
 ce qu'offre le th\'eor\`eme pr\'ec\'edent :
\begin{prop}
Soit $(X^{2n},\om, J, c)$ une vari\'et\'e symplectique r\'eelle de partie r\'eelle non vide, 
et $(s_{k}-_{k}$ une suite de sections AH, uniform\'ement transverse  et sym\'etrique. Alors il existe une constante
$C$ ind\'ependante de $k$, telle que le nombre de composantes de la partie r\'eelle  de $s_{k}^{-1}(0)$
 n'exc\`ede pas $C k^{\frac{n}{2}}$. 
\end{prop}
\bpr [ de la proposition 7] Soit $x$ un z\'ero quelconque de $s_{k}$. Puisque $s_{k}$ est 
transverse, il existe un $\ep >0$ tel qu'en $x$, $|\nabla s_{k}|\geq 2\ep\sqrt k$. Par ailleurs
$\nabla \nabla s_{k}$ est un $O(k)$, si bien que pour $\eta$ assez petit ind\'ependant de $k$, 
 la d\'eriv\'ee reste en norme sup\'erieure \`a $\ep \sqrt k$ sur la boule $B_{g_k}(x,\eta)$. 
Il s'ensuit (cf. [Au], partie 3.3) que sur cette boule, l'hypersurface est triviale. 
Par cons\'equent, chaque composante connexe de $s_{k}^{-1}(0)\cap \Rr X$ 
contient une boule de rayon $\eta/\sqrt k$, telle que deux de ces boules ne s'intersectent pas. Si $N$ est le nombre
de ces composantes, $N(\eta/\sqrt k)^n$ est donc major\'e par le volume de $\Rr X$, ce qui implique la contrainte 
 $N \leq C k^{\frac{n}{2}}$ pour une constante $C$ ind\'ependante de $k$.
\epr
\bpr[ du th\'eor\`eme 6] 
Nous avons d'abord besoin du lemme suivant :
\begin{lemme}
Soit $x\in \bl$. Alors il existe une application $\phi: B_g(x)\to \Cc^n$ v\'erifiant :

(i) $\phi(x) = 0$

(ii) $\phi^*\omega_0 = \omega$

(iii) $\phi(\bl) = \Rr ^n$

(iv) $\phi^*J_{0|\bl} = J_{|\bl}$

(v) $\phi$ se rel\`eve en un isomorphisme entre $L$ et le fibr\'e trivial sur $\Cc^n$. 

(vi) $\tilde c = \phi_*c$, o\`u $c_0$ est la conjugaison sur $\Cc^n$ et
v\'erifie  $d\tilde c = dc$ aux points de $\bl$.
\end{lemme}
 \bpr[ du lemme] Les cinq premières 
 propri\'et\'es  sont classiques (cf. [Au, Mu, Pr]). Quant \`a la derni\`ere, remarquons d'abord que  
$dc_{|T\bl} = d\tilde c_{|T\bl}. $ 
Ensuite,  si $\lambda$ est un vecteur
normal \`a $T\bl$, alors 
$J\lambda\in T\bl$, et donc 
$$d\tilde c (\lambda) = - d\tilde c(J^2 \lambda) = J d\tilde c (J\lambda) 
= J^2\lambda = -\lambda = dc (\lambda)$$ ce qui d\'emontre le \textit{(vi)}.
\epr
L'id\'ee est de prendre sur $\Rr X$ un r\'eseau dont la maille est de $g_{k}$-taille $D$, 
la constante $D$ \'etant ind\'ependante de $k$, puis
pour chaque sommet $x_{i}$, de construire une suite de sections $\tau_{i,k}$ suffisamment transverse, et dont
la partie r\'eelle du lieu d'annulation est localement et approximativement une hypersph\`ere 
incluse dans $ B_{g_k}(x_i,1)$.  
Ce qui suit traite en fait d'une situation plus g\'en\'erale. Soit $f : \Bb^{2n} \to \Cc$ une fonction
holomorphe, non singuli\`ere, v\'erifiant $\overline{f(\bar z) } = f(z)$ pour tout $z\in 2\Bb^{2n}$, 
et telle que la partie r\'eelle de son lieu d'annulation soit incluse dans la boule ouverte. L'exemple
dont on se servira est  
$$f(z_{1}, \cdots, z_{n})= z_{1}^2+\cdots + z_{n}^2-\frac{1}{2}.$$
 Dans ce mod\`ele local, 
la suite de sections $$\tau_{k}(z) = f(z\sqrt k)e^{-k|z|^2}$$ est  holomorphe et  sym\'etrique. Elle est
 de plus uniform\'ement $\ep $-transverse pour un $\ep $ donn\'e sur la boule. En effet, 
 remarquons d'abord le fait suivant : 
 $$\exists \eta>0 , \  |f|<\frac{\eta}{2eC} \implies  |df|>\eta, $$ 
 o\`u $C$ est elle que $|\nabla e^{-k|z|^2}|\leq C\sqrt k$. Si bien que
 si $|\tau_k|<\eta/2e^2C$, 
 alors 
$$|\na \tau_{k}| \geq  \sqrt k |df| e^{-k|z|^2} - |f||\nabla  e^{-k|z|^2}|\geq \sqrt k\eta e^{-1}/2,$$
ce qui d\'emontre la transversalit\'e uniforme de $\tau_k$. 
De plus, il est clair que la d\'eriv\'ee de $\tau_k$ est en norme inf\'erieure \`a $C'\sqrt k$, o\`u $C'$ ne d\'epend pas 
non plus de $k$.  Enfin, la restriction \`a $\Rr^n$
de la suite de sections est \'egalement transverse, pour les m\^emes raisons que pr\'ec\'edemment. 

Nous utilisons le lemme pr\'ec\'edent afin de rapatrier cette suite de sections sur notre vari\'et\'e symplectique. 
Nous noterons  $\tau_{i,k}$ la section $\phi_{i}^*\tau_{k}$, o\`u $\phi_{i}$ est 
le diff\'eomorphisme du lemme, centr\'e sur le point $x_{i}$. 
La suite de sections obtenue est AH, sym\'etrique, et $\ep$-transverse pour un $\ep$ qu'on peut choisir uniforme pour
tout point du r\'eseau puisque la partie r\'eelle $\Rr X$ est compacte.
Le proc\'ed\'e de cut-off traditionnel chez Donaldson ne vient pas perturber ces propri\'et\'es. En revanche, 
la transversalit\'e n'est plus vraie que sur une boule qu'on peut choisir de $g_{k}$-rayon $1$.    

Maintenant, nous faisons la somme de toutes ces sections. Si la maille $D$ du r\'eseau est 
suffisamment grande, mais ind\'ependante de $k$, la somme de toutes les autres contributions 
laisse transverse chacune des sections particuli\`eres. 
En effet, la norme $C^0$ (resp. $C^1$) de la somme de toutes les autres contributions 
que $\tau_{i,k}$ est major\'ee
par :
\begin{eqnarray*}
|\sum_{x_{j} \in \text{R\'eseau} \setminus \{x_{i}\}} \tau_{x_{j},k}| &\leq  &
C e^{-D} \\
|\sum_{x_{j} \in \text{R\'eseau} \setminus \{x_{i}\}} \nabla \tau_{x_{j},k}| &\leq  &
C e^{-D}\sqrt k. 
\end{eqnarray*}
Au total, pour $D$ assez grand (ind\'ependant de $k$) la transversalit\'e d'origine 
au point $x_{i}$ est pr\'eserv\'ee par les ajouts des autres
fonctions $\tau$ sur la boule $B_{g_{k}}(x_{i},1)$. La section continue donc \`a s'annuler sur $\Rr X$
sur une sous-vari\'et\'e isotope \`a la r\'eunion des hypersph\`eres, et ce dans un voisinage tubulaire
de celles-ci et de $g_{k}$-taille de l'ordre de $e^{-D}/\ep$, bien inf\'erieure \`a
 celle qui s\'epare $x$ des autres points du r\'eseau, si l'on choisit
$D$ assez grand. On a donc cr\'e\'e
autant d'hypersurfaces (d\'eform\'ees des hypersph\`eres) d'annulation que de points du r\'eseau, soit $\ep' \sqrt k$, o\`u 
$\ep'$ est une constante assez petite ne d\'ependant que de la g\'eom\'etrie de $(X,\om)$. Enfin, la proposition 2
nous permet de perturber $s_{k}$ en une suite de sections transverses et sym\'etriques. Si la perturbation
est assez faible, les pseudo-hypersph\`eres r\'eelles ne sont pas d\'etruites, et le th\'eor\`eme est d\'emontr\'e.
\epr
\noindent
\subsection{Le cas int\'egrable}
\bpr[ de la proposition  1]
La proposition 1 est d\'eduite de la conjonction des r\'esultats pr\'ec\'edents, et du fait fondamental
que lorsque $J$ est int\'egrable, c'est-\`a-dire que $X$ est une vari\'et\'e \kah, 
il est possible (cf. [Do1], proposition 34) de rendre les  sections
concentr\'ees du lemme 7 \textit{holomorphe}. Par ailleurs 
il est clair que si $s$ est une section holomorphe, $\ka(s)$ est \'egalement holomorphe. 
Dans notre cas, il est crucial de trouver, pour tout point $p$ de $\Rr X$, une section \textit{sym\'etrique} concentr\'ee
et holomorphe du type pr\'ed\'edent. Il suffit comme dans la partie 1.3 de changer la phase de $\si_{k,p}$, 
puis de prendre $\frac{1}{2}(\si_{k,p}+ \ka(\si_{k,p}))$, qui conserve toutes les propri\'et\'es souhait\'ees.
En utilisant la version holomorphe de ces sections concentr\'ees, on obtient tous les r\'esultats pr\'ec\'edents dans le cadre
complexe. 
\epr
\section{Pinceaux de Lefschetz r\'eels}
L'existence de pinceaux de Lefschetz r\'eels se fait en perturbant une paire de section AH. 
La proposition suivante (cf.
D\'efinition 5 de [Do2] et  Proposition 
5.5 de [Au,Mu,Pr]) repr\'esente la premi\`ere \'etape \`a r\'ealiser avant le th\'eor\`eme : 
\begin{prop} Soient $s_0$ et $s_1$ 
des suites de sections AH et sym\'etriques. Alors
pour $\ep>0$ assez petit, il existe une suite de  
sections $\tau_0\oplus \tau_1$ de $L^k\oplus L^k$ sym\'etriques, telle que :

1. La section $s_0 + \tau_0$ est $\ep$-transverse.

2. La section $(s_0+\tau_0)\oplus (s_1+\tau_1)$ est $\ep$-transverse,

3. La (1,0)-d\'eriv\'ee $\partial F$ de la fonction complexe $F =(s_1+\tau_1)/(s_0+\tau_0)$, 
est $\ep$-transverse sur l'ensemble $Z_{k,\ep} = \{ |s_0+\tau_0|\geq \ep\}$.

4. $F(c) = \bar F$. 
\end{prop}
\subsection{De la proposition 8 au th\'eor\`eme 4}
Dans ce paragraphe, nous supposons que la proposition pr\'ec\'edente est vraie,
et nous d\'emontrons sous cette hypoth\`ese que chacune des conditions
de la d\'efinition 1 est v\'erifi\'ee. 

\subsubsection{Propri\'et\'e (\textit{iv})}
Nous devons perturber les suites de sections $s_0$ et $s_1$ afin 
de satisfaire \`a la condition de la propri\'et\'e $(iv)$ de la d\'efinition 1. Cette perturbation
est r\'ealis\'ee dans [Do2], pp 214-215, Lemme 11. Comme d'habitude, 
nous devons nous assurer que nous pouvons la r\'ealiser de façon
sym\'etrique. Rappelons le principe de [Do2]. En un point $x$ de 
$N = \{s_0 = s_1 = 0\}$, l'espace $T_x N$ est le noyau de l'op\'erateur :
$$ D_x := \nabla s_0 \oplus \nabla s_1 : T_xX \to L^k_x\oplus L^k_x.$$
Soit $N_0$ le suppl\'ementaire symplectique de $T_xN$ dans $T_{x}X$. 
L'op\'erateur $D_x$ \'etablit un isomorphisme (r\'eel) entre
$N_0$ et $L_x^k \oplus L_x^k$, et induit ainsi par tir\'e en arri\`ere une 
structure complexe $j(D)$ sur $N_0$. 
\begin{lemme}(Lemme 11 de [Do2]) 
Pour tout point $x$ dans $N$, $F = s_1/s_0$ peut \^etre repr\'esent\'ee
sous la forme $(iv)$ de la d\'efinition 1 si et seulement si
la restriction de la forme symplectique $\om$
\`a $N_0$ est une forme (1,1) pour $j(D)$ et positive.
\end{lemme}
Nous d\'emontrons maintenant le lemme suivant :
\begin{prop}
Il existe une perturbation  de $s_0\oplus s_1$ 
de sorte que $F$ soit sym\'etrique et v\'erifie les conditions du lemme pr\'ec\'edent.
\end{prop}
\begin{preuve}
La d\'emonstration de la proposition se fait en deux temps. D'une part
il nous faut trouver une structure complexe $j$ sur chaque fibre $N_0$ sym\'etrique
par rapport \`a $c$, telle que $\om_{|N_0}$ soit $(1,1)$ pour
$j$. D'autre part, un th\'eor\`eme d'inversion locale 
nous permet de perturber sym\'etriquement les sections $s=s_{0}\oplus s_{1}$ en $\tilde s$ 
 de sorte que $j= j(D(\tilde s))$. 
Soit $$\pi : TX  \to N_0$$ la projection symplectique sur $N_0$. Par d\'efinition
de $N_{0}$, $g$, $\omega$ et $J$, l'endomorphisme $\pi J : N_0 \to N_0$ v\'erifie
$$g(u,v) = \om (u,\pi Jv) \  \forall u\in N_0, \ \forall v\in N_0.$$ 
Si $\pi J$ n'est pas tout \`a fait une structure complexe, ce n'est pas loin :
$$(\pi J)^2 = -Id_{|N_0} + O(1/\sqrt k). $$
En effet, $N_0$ et $TN$ sont approximativement des vari\'et\'es $J$-complexes, 
\`a $C/\sqrt k$ pr\`es, ce qui implique que  $$|J-\pi J|\leq C/\sqrt k,$$ et donc aussi 
 l'estimation pr\'ec\'edente. 
De plus, on peut perturber $a =\pi J$ en une authentique structure complexe
$j$ sur $N_0$, par la m\'ethode classique suivante. Il est clair que $a$ est 
antiautoadjoint relativement \`a  la m\'etrique $g_{|N_{0}}$. Soit $q$ l'endomorphisme
d\'efini comme la racine carr\'ee de $-aa^*$. L'application $j = aq^{-1}$
est alors une structure complexe compatible avec $\om_{|N_0}$.
 L'estimation pr\'ec\'edente montre que
$$|\pi J- j|\leq C/\sqrt k.$$
 Enfin,
chacun des \'el\'ements \'etant (anti)covariants par $dc$, le r\'esultat 
est \'egalement covariant par $dc$, i.e $c^*j =-j$. 

Montrons maintenant que cette structure est atteinte 
par l'interm\'ediaire d'une perturbation des sections de d\'epart. 
Pour cela, soit $\bj_0$ l'ensemble des structures complexes sur $N_0$.  l'application :
\begin{eqnarray*}
j : Gl (N_0, L^k\oplus L^k) & \to          &\bj_0         \\ 
                       D                  & \mapsto  & D^{-1}i D,
\end{eqnarray*}
o\`u $i$ est la structure complexe sur $L^k\oplus L^k$.
 Puisque $s_0\oplus s_1$ est asymptotiquement 
holomorphe, on a $$||j(\nabla s))-\pi J||  \leq C |(\nabla s)^{-1}|\leq C/(\eta \sqrt k),$$
puisque  la norme de l'inverse de $D_x$ est inf\'erieure \`a $1/(\eta\sqrt k)$, en vertu de la $\eta-$transversalit\'e
de la section $s =s_0\oplus s_1$, et que $s$ est AH, donc son $\del$ est inf\'erieur en norme \`a $C$. 
La diff\'erentielle de $j$ en $D$ est 
$$d_Dj(H) = -D^{-1}HD^{-1}iD + D^{-1}iH = 2D^{-1}H^{0,1}j(D),$$
o\`u la partie $(0,1)$ de l'op\'erateur $H$ est entendue en fonction des  structures
complexes $j(D)$ et $i$. Cette diff\'erentielle est surjective,
car $$T_{j(D)} \bj_{0} = \{K, \ Kj(D) + j(D) K = 0\},$$
et $K=d_Dj(H)$ \'equivaut \`a $H^{0,1} = -\frac{1}{2}DKj(D)$, \'equation qui poss\`ede une solution car 
 l'op\'erateur $DKj(D)$ est $(0,1)$.
De plus, il existe des constantes $\epsilon >0$ et $c$ ind\'ependantes de $k$, telles que 
$$ \nu(d_{D(s)}j) \geq \epsilon \ \text{et} \ |j|_{C^2}\leq c.$$

Enfin, on a vu pr\'ec\'edemment que  $|j(D(s))-j|\leq C/\sqrt k$. 
En conclusion, le th\'eor\`eme des fonctions
implicites nous donne l'existence d'un op\'erateur $\tilde D$ tel que 
$$j(\tilde D) = j \ \text{et}\ |\tilde D-D(s)| \leq C/\sqrt k.$$
Dans notre cas, au lieu de consid\'erer les sections du fibr\'e $Gl(N_0, L^k\oplus L^k)$, on
choisit les sections $dc$-\'equivariantes de ce fibr\'e, c'est-\`a-dire les op\'erateurs $D$ tels que
$$D_{c(x)}d_xc = d_{c(x)}D_x.$$ L'application $j$ \'etant \'equivariante, 
le r\'esultat est qu'il est possible de trouver un $D$ \'equivariant r\'esolvant 
$j(D) = j$. 
Maintenant, il suffit de perturber $s = s_0\oplus s_1$ en $\tilde s$ de façon \'equivariante
sur un voisinage de $N$ de taille constante (pour
la m\'etrique $g_k$) afin que le nouveau $D(\tilde s) $ soit \'egal \`a $\tilde D$. Ceci
peut se faire sans changer l'ensemble $N$.
\end{preuve}

\subsubsection{Propri\'et\'e $(v)$}
Nous suivons maintenant les pages 209 \`a 213 de [Do2], dont le contenu
permet de perturber la fonction $F=s_1/s_0$ de sorte qu'elle convienne
au mod\`ele $(v)$ de la d\'efinition 1. Pour cela, rappelons le contenu de la proposition
9 de [Do2] :
\begin{prop} Soit $\Delta = \{\partial F = 0\}$, $\Gamma=\{|\partial F| \leq |\del F|\}$,
et enfin $\Omega_\chi = \{|s_0|\geq \chi\}$. Alors 
\begin{enumerate}
\item  $\Delta$ est un ensemble fini.
\item Il existe $\chi$ et $\rho_0$ ind\'ependants de $k$ tel que  
les $g_k$-boules centr\'ees sur les points de $\Delta$ et de taille $\rho_{0}$
soient disjointes et contenues dans $\Omega_\chi$.
\end{enumerate}
\end{prop}
Fixons  $x$ un point de $\Delta$. 
La d\'eriv\'ee covariante (induite par la connexion
de Levi-Civita) de la 1-forme $\partial F$ se d\'ecompose de la façon suivante :
$$\nabla (\partial F) = \partial_\nabla (\partial F) + \bar \partial_\nabla (\partial F).$$
Puisque l'involution est anti-holomorphe et isom\'etrique, et que $F(c) = \bar F$, 
on obtient facilement que 
$$\overline{c^*\partial_\nabla (\partial F)} = \partial_\nabla (\partial F).$$ 
Soit $H(z) = \sum H_{\al \beta}z_\al z_\beta$ un polyn\^ome complexe
en des coordonn\'es adapt\'ees au point $x$ \`a $\om$, et tel que
sa hessienne en $x$ soit \'egale \`a $\partial_\nabla (\partial F)$. 
La perturbation de $F$ dans le cas classique  sur $B_{g_{k}}(x,\rho)$ se fait de la façon suivante : 
$$ \tilde F(z) =\beta_\rho (w+H(z)) + (1-\beta_\rho)F(z),$$
o\`u $\beta_{\rho}$ est une fonction plateau \`a support dans $B_{g_{k}}(x,\rho)$
Le lemme de [Do2] d\'ecrit l'effet de cette modification :
\begin{lemme} Pour $\rho$ assez petit, $k=k(\rho)$
assez grand, et $w$ assez proche (relativement \`a $\rho$) de $F(x)$,
alors $x$ est le seul point de $B_{g_{k}}(x,\rho)$ o\`u $|\partial F|\leq |\del F|$.
\end{lemme}

Si $x$ appartient \`a $\Delta \cap \Rr X$, 
$F(x)$ est r\'eel, si bien qu'il suffit de choisir $w$ r\'eel pour que la perturbation pr\'ec\'edente r\'ealise la condition $(v)$,
tout en laissant $F$ sym\'etrique. Si $x$ appartient \`a  $\Delta$ mais pas \`a $\Rr X$, 
on transforme 
$F$ en $$ \tilde F(z) =\beta_\rho (w+H(z)) +
 (1-\beta_\rho)F(z) + c^*(\beta_\rho (w+H(z)) + (1-\beta_\rho)F(z)).$$
La proposition 10 montre que si $x$ et $c(x)$ sont diff\'erents, alors ils sont forc\'ement 
$g_k$-distants d'au moins $\rho_0$. Le  
lemme pr\'ec\'edent montre alors que si $\rho$ est assez petit, les supports des deux
 parties de la fonction perturbatrice ne se rencontrent pas, si bien qu'on obtient
 sans effort la propri\'et\'e $(v)$ en $x$ et $c(x)$. 
Par ailleurs en choisissant $w$ de sorte que $F(x)+ w$ ne soit pas r\'eel, 
on a $\tilde F(x)\not= \tilde F(c(x))$. 

\subsection{D\'emonstration de la proposition 8}
\noindent
\textbf{Les propri\'et\'es 1 et 2.}
La propri\'et\'e 1 a \'et\'e r\'ealis\'ee dans la premi\`ere partie. 
Pour la propri\'et\'e 2., nous reprenons les arguments de la partie 3.3 dans [Au], 
(cf. aussi [Au,Mu,Pr], ainsi que notre proposition 2 avec $E=\Cc^2$). Le lieu des z\'eros $Z_{0,k}$ de $s_0$ 
est une hypersurface symplectique r\'eelle, 
telle que la distance entre $TZ$ et $JTZ$ 
est en $1/\sqrt k$.  On peut alors
appliquer la premi\`ere partie \`a la suite de sections
$s_{1|Z_{0,k}}$, c'est-\`a-dire la rendre 
transverse sur $Z_{0,k}$ tout en la laissant
sym\'etrique pour la structure r\'eelle induite par $c$.
La morphologie de $Z_{0,k}$ varie avec $k$, mais cela ne pose
pas de probl\`emes (cf. [Au] pour les d\'etails). 
Maintenant, si la restriction de $s_1$ est transverse, 
$s_1$ est \'egalement transverse dans l'espace ambiant aux points
de $Z_{0,k}$. Enfin, on peut constater qu'alors la section 
perturb\'ee $s_0\oplus s_1$ est transverse sur un voisinage fixe
de $Z_{0,k}$.\\
 
\noindent
\textbf{Propri\'et\'e 3 : la transversalisation de $\partial (s_1/s_0)$.}
Il reste \`a transversaliser la section $f=\partial (s_1/s_0)$. 
Cela se fait sur trois types d'endroits : hors d'un voisinage
de $g$-taille fixe, sur un $g_k$-voisinage
fixe de $\Rr X$, et enfin entre les deux pr\'ec\'edents.
Dans le premier cas, nous suivons la construction de Donaldson. 
Etant donn\'e que d'une part  les supports des sections $s_k$ et $\kappa (s_k)$ sont
disjoints, et que d'autre part la condition de transversalit\'e 
de $\partial (s_1/s_0)$ est invariante par $c$, il suffit
de r\'ealiser la perturbation sur la moiti\'e des boules, 
et de perturber par $g+\overline {c^*g}$ plut\^ot que simplement par $g$. 
La transversalisation dans la derni\`ere situation interm\'ediaire est r\'ealis\'ee dans [Au,Mu,Pr] (p. 22). 
Leur travail s'adapte imm\'ediatement \`a notre cas. 
En ce qui concerne la seconde situation, dans le cas trait\'e par 
[Au,Mu,Pr], la paire de sections est telle que la restriction de leur rapport 
est une fonction de Morse donn\'ee sur $\Rr X$, si bien que la transversalisation 
de $f$ est d\'ej\`a faite. 
Dans notre cas, nous partons de deux sections quelconques, si 
bien qu'il est n\'ecessaire de transversaliser pr\`es de $\Rr X$. N\'eanmoins, 
nous utilisons essentiellement la m\'ethode des trois auteurs pr\'ec\'edemment cit\'es. 
La transversalisation de $f$ par  [Au,Mu,Pr] utilise une r\'ecurrence de la façon suivante. Soit $U$ un ouvert connexe,
simplement
connexe et invariant par $c$. 
Choisissons des sections $\alpha_{i}$ 
orthonormales du fibr\'e trivial
$T^{1,0}X_{|U}$, et nommons $A_r$ le fibr\'e en droites complexes engendr\'ees par les valeurs de $\alpha_{r}$. On a donc
$\bigoplus_{r=1}^{r=n} A_r= T^{1,0}X_{|U}$. Soit de plus $\pi_{\leq r}$ la projection
de $T^{1,0}X$ sur $ E_r = \bigoplus_{i\leq r} A_r$, et $\pi_{r}$ la projection sur $A_{r}$. 

Le principe de r\'ecurrence utilise le lemme
suivant :
\begin{lemme}([Au]) Si $s$ est une section de $E\otimes L^k$ transverse sur $W=\{s=0\}$,
et $\tau$ une section AH de $L^k$ telle que  
$t_{|W}$ est transverse sur $W$, alors $s\oplus t$ est transverse.
\end{lemme}
Ainsi, on suppose que $\pi_{\leq r} f$ est transverse sur $W_{r} = \{\pi_{\leq r} f= 0\}$. 
On tente ensuite de perturber $F$ pour obtenir la transversalisation
sur $W_r$ de $\pi_{r+1} f$. Au bout de $n$ it\'erations, on obtient
la transversalisation de $F$ sur l'ouvert. 
Supposons donc que $\pi_{\leq r} f$ soit transverse sur $W_r$.
Nous recouvrons l'intersection de $W_{r}$ avec le voisinage
de $\Rr X$  par des boules centr\'ees sur des points de $\Rr X$ (et non pas
sur des points de $W_{r}$ comme dans [Au,Mu,Pr]).

 Soit $x\in \Rr X$ un de ces centres, et $(z_i)$ des coordonn\'ees complexes
 centr\'ees $x$, telles que pour tout $i$, $ \partial z_i(x) = \alpha_i (x)$.
Puis perturbons $s_1$ de la façon suivante :
$$ s'_1 = s_1 + wz_{r+1} \siref.$$
On a alors 
$\partial (s'_1/s_0) = \partial (s_1/s_0) + 
w (\siref /s_0\partial (z_{r+1})  + z_{r+1} \partial (\siref /{s_0})).$
 Rappelons qu'il est n\'ecessaire de transversaliser $f$ uniquement sur
  $X_{\ep, k}=\{|s_{0}|\geq \ep\}$. 
  Si on appelle $\zeta$ la section de $A_{r+1}$ d\'efinie par 
 $$\zeta = (\siref/s_{0}) \pi_{r+1}(\partial z_{r+1})
                                          +  z_{r+1}\pi_{r+1}(\partial (\siref/s_{0}))),$$
                          alors $\zeta $ est uniform\'ement minor\'ee
                  sur $W_{r}\cap X_{\ep, k}\cap B_{g_k}(x,\rho)$, 
       si $\rho$ est suffisamment petit (par rapport \`a $\ep$ et ind\'ependamment de $k$ et $x$). 
      En effet, l'estimation       $ |\partial (\siref/s_{0}) |\leq C\sqrt k$ ($C$ d\'epend
      cette fois de $\ep$)
                   montre que pour $\rho$ assez petit (par rapport \`a $\ep$),
                   la seconde moiti\'e de $\zeta$ est inf\'erieure en norme \`a 1/10
                   de  la norme de la premi\`ere partie, 
                   qui est minor\'ee par une constante strictement positive d\'ependant de $\ep$.                                         
                    La section 
 $$
\pi_{r+1} (\partial (s_{1}'/s_{0}))= \pi_{r+1}(\partial (s_{1}/s_{0})) + w\zeta
 $$ 
 est donc bien d\'efinie sur l'intersection de  $W_{r}$ avec la boule $ B_{g_k}(x,\rho) $ et $X_{\ep, k}$,
 ainsi que la fonction
$$f_{w} = \frac{\pi_{r+1}\partial (s_1/s_0)}{\zeta} +w. $$
Pour toute composante connexe $C$ de l'intersection de $W_{r,w}$ avec
la boule $B_{k}(x)$, on choisit une trivialisation (cf. [Au] et [Au, Mu, Pr]) 
de $C$. La restriction \`a celle-ci de la fonction poss\`ede les estim\'ees d'holomorphie asymptotique. 
Maintenant la proposition 4 permet de trouver un $w$ r\'eel tel que
$f_{w| C }$ devienne $\eta$-transverse sur $C$. On recommence le processus pour toutes les composantes,
en nombre born\'e ind\'ependant de $k$ (mais d\'ependant de $\ep$), pour transversaliser $f_{w|W_{r,w}}$. 
Au total, la section $$\pi_{r+1}(\partial (s_1/s_0)+w \pi_{r+1}\partial (z_{r+1} \siref/s_0)$$
est sym\'etrique sur $X$, et transverse sur $W_{r,w}$. Par r\'ecurrence,
on obtient donc le r\'esultat.


\begin{thebibliography}{papousxxx}
\bibitem[Au]{Do}
    D.  \textsc{Auroux},
     \emph{Asymptotically holomorphic families of symplectic submanifolds}, 
GAFA. 7 (1997) 971-995.
\bibitem[Au,Ga,Mo]{Do}
    D.  \textsc{Auroux}, D. \textsc{Gayet}, J.-P \textsc{Mohsen},
     \emph{Symplectic hypersurfaces in the complement of an
isotropic submanifold}, Math. Ann. 321, 739-754 (2001).
\bibitem[Au,Mu,Pr]{Do}
    D.  \textsc{Auroux}, V. \textsc{Mu\~{n}oz}, F. \textsc{Prezas},
     \emph{Lagrangian submanifolds and Lefschetz pencils}, Arxiv:math.SG/0407126.
\bibitem[Do1]{Do}
    S.K.  \textsc{Donaldson},
     \emph{Symplectic submanifolds and almost-complex geometry}, J. 
Diff. Geom.Vol. 44 (1996),
666-705.
\bibitem[Do2]{Do}
    S.K.  \textsc{Donaldson},
     \emph{Lefschetz pencils on symplectic manifolds}, J. 
Diff. Geom.Vol. 53 (1999),
205-236.


\bibitem[Ed, Ko]{Do}
    D.  \textsc{A. Edelman, E. Kostlana}, 
      \emph{How many zeros of a random polynomial are real ?},
      Bull. of the A.M.S, Vol. 32 (1995), pp 1-37.  
 
\bibitem[We]{We}
     J-Y. \textsc{Welschinger},
     \emph{Invariants of real symplectic 4-manifolds and lower bounds 
in real enumerative geometry}, Inv. Math. Vol. 162 (2005).

\bibitem[We2] {We}J-Y. \textsc{Welschinger},
     \emph{Real structures on minimal ruled surfaces}, 
Comment. Math. Helv. 78, No. 2, 418-446 (2003). 
 \end{thebibliography}
\end{document}